\documentclass[a4paper,12pt]{article}
\usepackage[utf8]{inputenc}
\usepackage{amssymb,amsmath,amsthm}
\usepackage{hyperref}
\usepackage{tikz}
\usepackage{graphicx}
\newtheorem{theorem}{Theorem}

\newtheorem{lemma}[theorem]{Lemma}

\newtheorem{conjecture}[theorem]{Conjecture}

\newcommand{\ceil}[1]{\left\lceil{#1}\right\rceil}

\usepackage{authblk}
\DeclareMathOperator{\ex}{ex}

\title{A new construction for the planar Tur\'an number of cycles}
\author[1]{Ervin Gy\H{o}ri}
\author[1,2,3]{Kitti Varga}
\author[1,4]{Xiutao Zhu}

\affil[1]{Alfr\'ed R\'enyi Institute of Mathematics, Hungarian Academy of Sciences.}
\affil[2]{Department of Computer Science and Information Theory, Budapest University of Technology and Economics.}
\affil[3]{ELKH-ELTE Egerv\'ary Research Group.}
\affil[4]{Department of Mathematics, Nanjing University.}

\date{}

\begin{document}

\maketitle

\begin{abstract}
The planar Tur\'an number $\ex_{\mathcal{P}}(n,C_k)$ is the maximum number of edges in an $n$-vertex planar graph not containing a cycle of length~$k$. Let $k\ge 11$ and $c,d$ be constants. Cranston et al., and independently Lan and Song showed that $\ex_{\mathcal{P}}(n,C_k)\ge 3n-6- cn/k$ holds for large $n$. Moreover, Cranston et al.\ conjectured that $\ex_{\mathcal{P}}(n,C_k)\le 3n-6- dn/k^{\log_2 3}$ when $n$ is large.

In this note, we prove that $\ex_{\mathcal{P}}(n,C_k)\ge 3n-6-6\cdot 3^{\log _23}n/k^{\log_2 3}$ holds for every $k\ge 7$. This implies that Cranston et al.'s conjecture is essentially best possible.

\end{abstract}

\section{Introduction}

In this article, the cycle and complete graph on $k$ vertices are denoted by $C_k$ and $K_k$, respectively. For $k\in \{4,5\}$, let $\Theta_k$ denote the graph obtained from $C_k$ by adding a chord.

The Tur\'an number $\ex(n, H)$ for a graph $H$ is the maximum number of edges in an $n$-vertex graph containing no copy of $H$ as a subgraph. The first result on this topic was obtained by Tur\'an~\cite{Turan}, who proved that the balanced complete $r$-partite graph is the unique extremal graph for $\ex(n,K_{r+1})$. The Erd\H{o}s--Stone--Simonovits Theorem~\cite{Erdos1966ALT,erdos1946structure} generalizes this result  and  asymptotically  determines  $\ex(n,H)$  for  all  non-bipartite graphs~$H$:  $\ex(n,H)= \big(1-\frac{1}{\chi(H)-1} \big)\binom{n}{2}+o(n^2)$, where $\chi(H)$ denotes the chromatic number of $H$.

In 2016, Dowden~\cite{Dowden} initiated the study of planar Tur\'an problems, in which we want to determine the maximum number of edges of an $n$-vertex planar graph containing no copy of  $H$ as a subgraph. This number is denoted by $\ex_{\mathcal{P}}(n,H)$.
Dowden proved that $\ex_{\mathcal{P}}(n,C_4)\le (15n-30)/7$ for all $n\ge 4$, and $\ex_{\mathcal{P}}(n,C_5)\le (12n-33)/5$ for all $n\ge 11$~\cite{Dowden}.   Lan, Shi and Song~\cite{2017Extremal} showed that $\ex_{\mathcal{P}}(n,\Theta_4)\le (12n-24)/5$ for all $n\ge 4$, and $\ex_{\mathcal{P}}(n,\Theta_5)\le (5n-10)/2$ for all $n\ge 5$, and $\ex_{\mathcal{P}}(n,C_6)\le (18n-36)/7$. While the bounds on $\ex_{\mathcal{P}}(n,\Theta_4)$ and $\ex_{\mathcal{P}}(n,\Theta_5)$ are tight for infinitely many $n$. The upper bound on $\ex_{\mathcal{P}}(n,C_6)$ was improved by Ghosh, Gy\H{o}ri, Martin, Paulos and Xiao~\cite{6cycle}. They proved $\ex_{\mathcal{P}}(n,C_6)\le (5n-14)/2$ for all $n\ge 18$, and they also proposed the following conjecture.

\begin{conjecture}[Ghosh et al.~\cite{6cycle}]\label{Conj1}
For each $k\ge 7$ and sufficiently large $n$, 
\[\ex_{\mathcal{P}}(n,C_k)\le 3n-6-\frac{3n+6}{k}.\]
\end{conjecture}

Recently,  Conjecture 1 was disproved by Cranston, Lidick\'y, Liu and Shantanam~\cite{2021Planar} and independently by Lan and Song~\cite{LanSong} for $k\ge 11$ and sufficiently large $n$. (Cranston et al.'s lower bound is a little weaker than Lan and Song's)
\begin{theorem}[Lan and Song~\cite{LanSong}]
 Let $k\ge 11$ and $n\ge k-4+\lfloor(k-1)/2\rfloor$. Then there is a constant $c_k$ such that 
 \[\ex_{\mathcal{P}}(n,C_k)\ge \left(3-\frac{3-\frac{2}{k-2}}{k-6+\lfloor (k-1)/2\rfloor}\right)n+c_k.\]
\end{theorem}

Furthermore, Cranston et al.\ proposed a revised conjecture.
\begin{conjecture}[Cranston et al.~\cite{2021Planar}]\label{Conj2}
There exists a constant $d$ such that for all $k$ and all sufficiently large $n$, we have
\[\ex_{\mathcal{P}}(n,C_k)\le 3n-6-\frac{dn}{k^{\log_23}}.\]
\end{conjecture}

In this note, we give a new construction for the lower bound of $\ex_{\mathcal{P}}(n,C_k)$ and obtain the following theorem.

\begin{theorem}\label{thm2}
 For all $k\ge 7$ and sufficiently large $n$, we have 
\[\ex_{\mathcal{P}}(n,C_k)\ge 3n-6-\frac{6\cdot 3^{\log_23}n}{k^{\log_23}}.\]
\end{theorem}

Note that Theorem~\ref{thm2} implies Conjecture~\ref{Conj2} is essentially best possible.

\section{Our construction}

In this section, we show our construction and prove Theorem \ref{thm2}. We first define a sequence of planar graphs $T_i$ and use them in our construction.  This sequence was first introduced by Moon and Moser~\cite{1963Simple}.
 
Let $T_1$ be a copy of $K_4$ and $xyz$ be the outer cycle. Suppose $T_{i-1}$ is defined for some $i\ge 2$. Let $T_i$ be the graph obtained from $T_{i-1}$ as follows: in each inner face of $T_{i-1}$, add a new vertex and join the new vertex to the three vertices incident to this face.

By the above construction, each $T_i$ is a triangulation, i.e., a planar graph whose each face is a triangle. The outer cycle of $T_i$ is $xyz$ and
\[|V(T_i)|=4+3+3^2+\cdots+3^{i-1}=\frac{3^i+5}{2}.\]
Furthermore, Chen and Yu \cite{ChenGuantao2002Long} showed that $T_i$ has the following properties. 

\begin{lemma}[Chen--Yu~\cite{ChenGuantao2002Long}]\label{Chen-Yu}
 For any integer $i \ge 2$, we have the following.
 
 \begin{enumerate}
  \item[(i)] The length of the longest path between $x$ and $y$ in $T_i$ is $3\cdot 2^{i-1}$.
  \item[(ii)] The length of the longest cycle in $T_i$ is $7\cdot 2^{i-2}$.
\end{enumerate}   
\end{lemma}

\begin{figure}[ht]
 \begin{center}
 \begin{tikzpicture}[scale=1.25]
  \tikzstyle{vertex}=[draw,circle,fill=black,minimum size=5,inner sep=0]
  
  \coordinate (u) at (0,2);
  \coordinate (v) at (0,-2);
  \coordinate (w1) at (-3,0);
  \coordinate (w2) at (-2,0);
  \coordinate (w3) at (-0.5,0);
  \coordinate (w4) at (0.5,0);
  \coordinate (w5) at (2,0);
  \coordinate (w6) at (3,0);
  
  \fill[black!30, opacity=0.5] (u) -- (w1) -- (v) -- (w2) -- (u);
  \fill[black!30, opacity=0.5] (u) -- (w3) -- (v) -- (w4) -- (u);
  \fill[black!30, opacity=0.5] (u) -- (w5) -- (v) -- (w6) -- (u);
  
  \node at (-2.375,0) {$H_1^-$};
  \node at (0,0) {$H_j^-$};
  \node at (2.375,0) {$H_s^-$};
  
  \node[vertex] at (u) [label={above:$x$}] {};
  \node[vertex] at (v) [label={below:$y$}] {};
  \node[vertex] at (w1) [label={[xshift=-9pt, yshift=-12pt] $w_1$}] {};
  \node[vertex] at (w2) [label={[xshift=11pt, yshift=-12pt] $z_1$}] {};
  \node[vertex] at (w3) [label={[xshift=-10pt, yshift=-13.5pt] $w_j$}] {};
  \node[vertex] at (w4) [label={[xshift=9pt, yshift=-13.5pt] $z_j$}] {};
  \node[vertex] at (w5) [label={[xshift=-11pt, yshift=-12pt] $w_s$}] {};
  \node[vertex] at (w6) [label={[xshift=9pt, yshift=-12pt] $z_s$}] {};
  
  \draw[fill] (-0.725,1) circle (0.3pt);
  \draw[fill] (-0.625,1) circle (0.3pt);
  \draw[fill] (-0.525,1) circle (0.3pt);
  
  \draw[fill] (-0.725,-1) circle (0.3pt);
  \draw[fill] (-0.625,-1) circle (0.3pt);
  \draw[fill] (-0.525,-1) circle (0.3pt);
  
  \draw[fill] (0.525,1) circle (0.3pt);
  \draw[fill] (0.625,1) circle (0.3pt);
  \draw[fill] (0.725,1) circle (0.3pt);
  
  \draw[fill] (0.525,-1) circle (0.3pt);
  \draw[fill] (0.625,-1) circle (0.3pt);
  \draw[fill] (0.725,-1) circle (0.3pt);
  
  \draw[thick] (u) -- (w1) -- (v);
  \draw[thick] (u) -- (w2) -- (v);
  \draw[thick] (u) -- (w3) -- (v);
  \draw[thick] (u) -- (w4) -- (v);
  \draw[thick] (u) -- (w5) -- (v);
  \draw[thick] (u) -- (w6) -- (v);
  
  \draw[thick] (u) to [bend right=90, looseness=3.25] (v);
 \end{tikzpicture}
 \caption{The graph $H$.} \label{figure}
 \end{center}
 \end{figure}
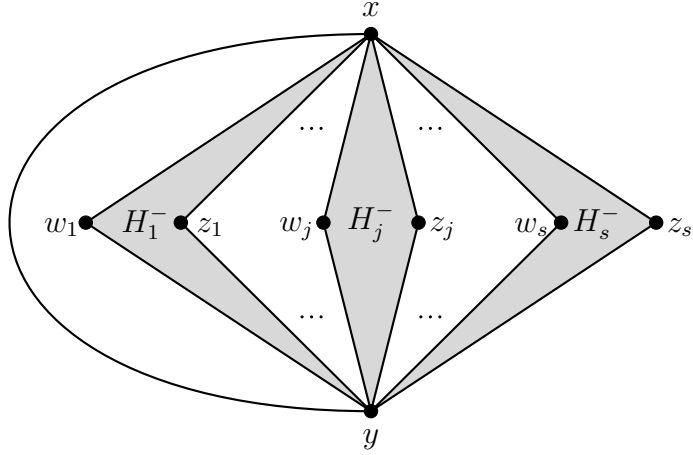

\vskip 2mm
After defining the graphs $T_i$, we show our construction.  Let $k\ge 7$ and $i$ be the maximum integer such that
\begin{align}\label{eq1}
3\cdot 2^{i-1}< \frac{k}{2}, \text{~~i.e.,~~} i=\big\lceil\log_2k/3\big\rceil-1.   
\end{align}
Suppose $n\ge \frac{3^{i}+5}{2}$ and $s= \big\lceil \frac{n-2}{(3^i+5)/2-2} \big\rceil$. 
Let $H_1,\ldots,H_{s-1}$ be $s-1$ copies of $T_i$ and $H_{s}$ be a subgraph of $T_i$ such that $H_s$ is a triangulation and has $n-(s-1)(\frac{3^i+5}{2}-2)$ vertices. Such a graph $H_s$ exists because of the process of the construction of $T_i$. For each $1\le j\le s$, we may assume $x_jy_jz_j$ is the outer cycle of $H_j$, and let the other triangular face containing $x_jy_j$ be $x_jy_jw_j$. Let $H_j^-$ be the subgraph obtained from $H_j$ by deleting the edge $x_jy_j$. Then each face of $H_j^-$ is a triangle except the outer face whose boundary is the cycle $x_jw_jy_jz_j$. Let us identify the vertices $x_j$ of $H_j^-$ for all $1 \le j \le s$ as a new vertex $x$, and identify the vertices $y_j$ of $H_j^-$ for all $1 \le j \le s$ as a new vertex $y$, and  add a new edge $xy$ to this graph. Let $H$ denote the resulting graph (see Figure~\ref{figure}).

We call each $H_j^-+xy$ a block of $H$. Note that each block is a copy of $T_i$ except the last block $H_{s}^-+xy$, which is a subgraph of $T_i$. Now we show that $H$ contains no~$C_k$. Let $C$ be a longest cycle in $H$. Since $\{x,y\}$ is a vertex-cut of $H$, $C$ passes through at most two blocks. If $C$ passes through only one block, then by (ii) of Lemma~\ref{Chen-Yu} and by the equation \ref{eq1}, we have 
\[ |V(C)| \le 7\cdot2^{i-2}<\frac{7k}{12}. \]
If $C$ passes through two blocks, then clearly $x,y\in V(C)$, thus by (i) of Lemma~\ref{Chen-Yu} and the equation \ref{eq1}, we have
\[ |V(C)| \le 2(3\cdot 2^{i-1})<k. \]
Hence, $H$ contains no $C_k$.

Next we calculate the number of edges of $H$. Clearly, $H$ is an $n$-vertex planar graph.
 By adding the edges $w_1z_2,~w_2z_3,\ldots, w_{s-1}z_{s}$ to $H$, it becomes a triangulation.  Hence, we have
 \[|E(H)|=3n-6-(s-1). \]
On the other hand, since $s= \big\lceil \frac{n-2}{(3^i+5)/2-2} \big\rceil$ and by equation \ref{eq1}, we can deduce 
\[\begin{split}
 |E(H)|\ge & 3n-6-\frac{2(n-2)}{3^i+1}=3n-6-\frac{2(n-2)}{3^{\ceil{\log_2k/3}-1}+1}\\
 \ge & 3n-6-\frac{6(n-2)}{3^{\log_2k/3}+3}\ge 3n-6-\frac{6\cdot 3^{\log_23}n}{k^{\log_23}}.
\end{split}\]
Hence, $\ex_\mathcal{P}(n,C_k)\ge 3n-6-\frac{6\cdot 3^{\log_23}n}{k^{\log_23}}$ and we are done.
 $\hfill\blacksquare$

\bibliography{citation.bib}

\end{document}